\documentclass[a4paper, 11 pt, twoside]{amsart}
\usepackage{srollens-en, shortcuts-sr}
\usepackage{todonotes}
\usepackage[left = 3.5cm, right = 3.5cm, headsep = 6mm,
footskip = 10mm, top = 35mm, bottom = 35mm, footnotesep=5mm, headheight =
2cm]{geometry}

\newlist{prooflist}{description}{1}
\setlist[prooflist]{font=\normalfont \itshape, labelindent = \parindent, leftmargin = 0pt}

\usepackage{comment}

\usepackage[unicode,bookmarks, pdftex]{hyperref}
\hypersetup{colorlinks=true,citecolor=NavyBlue,linkcolor=NavyBlue,urlcolor=Orange, pdfpagemode=UseNone, breaklinks=true}

\DeclareMathOperator{\Res}{Res}
\newcommand{\Nklt}{\mathsf{Nklt}}

\title{\texorpdfstring{$5I(K_X+\Delta)$}{5I(K+Delta)} is very ample on stable log surfaces}

\author{Wenfei Liu}
\address{School of Mathematical Sciences, Xiamen University, Siming South Road 422, 361005 Xiamen, Fujian Province, P.~R.~China}
\email{wliu@xmu.edu.cn}

\author{S\"onke Rollenske}
\address{S\"onke Rollenske
\\
FB 12/Mathematik und Informatik\\
Philipps-Universit\"at Marburg\\Hans-Meerwein-Str. 6\\
35032 Marburg\\Germany}

\email{rollenske@mathematik.uni-marburg.de}

\begin{document}
\begin{abstract}
We prove that Bombieri's bound for canonical models of surfaces of general type extends to stable log surfaces, that is,  $mI(K_X+\Delta)$ is very ample for every $m\geq 5$ and every stable log surface $(X,\Delta)$ over $\IC$ with reduced boundary and global Cartier index $I$. 

On the way to this result we determine the cases where $|3I(K_X+\Delta)|$ has base points. 
\end{abstract}

\maketitle
\setcounter{tocdepth}{1}
\tableofcontents

\section{Introduction}
Modular compactifications of moduli spaces are a central topic in algebraic geometry and for canonically polarised varieties the so-called KSBA-compactification parametrising stable varieties has been constructed \cite{KollarModuli}. 
In the case of surfaces it compactifies the classical Gieseker moduli space of canonical models \cite{gieseker77}. 

A crucial ingredient there, and in general an invaluable tool in the study of surfaces of general type, is Bombieri's study of pluricanonical maps \cite{bombieri73} and their later refinements (see e.g.~the survey \cite{BCP06}). In \cite{LR14}, the authors extended these methods
to stable log surfaces with reduced boundary (see Section \ref{sec: stable log surface} for definitions). At the time, the result seemed roughly as expected---the bounds a little worse than the classical ones as it happens for curves\footnote{For example, if a stable curve has an elliptic tail, then the canonical map is not a morphism.}---but we were unable to construct examples to show that the obtained  bounds were sharp.

To our surprise, Bombieri's bound extends to all stable log surfaces.
\begin{thm}\label{thm: A}
	Let $(X, \Delta)$ be a stable log surface over $\IC$ with reduced boundary $\Delta$ and  global Cartier index $I$. Then $mI(K_X+\Delta)$ is very ample for every integer $m\geq5$. 
\end{thm}
The proof follows a similar pattern compared to \cite{LR14}, combining Reider type results obtained by Kawachi, Ma\c{s}ek, and Langer \cite{Kawachi98,Masek99, KM98, Langer01} at normal points,  the classical method of restrictions to suitable pluricanonical curves and Fujino's method of  restriction to log-canonical centres 
\cite{Fujino14, Fuj17}. 
 The proof is further simplified by reducing the case with boundary to the case $\Delta=0$ via doubling along the boundary (Lemma~\ref{lem: double}). 
 
A crucial new ingredient is a characterisation of the base-point-freeness of the tricanonical map, made possible by the classification of Gorenstein stable (log) surfaces of volume one in \cite{FPR15, FPR15b, FPR17}.
\begin{thm}\label{thm: B}
Let $(X,\Delta)$ be a stable log surface with reduced boundary and global
index $I$. Then $|3I(K_X+\Delta)|$ has no base points unless $\Delta=0$, $I = K_X^2 = 1$ and $\chi(X)\in\{1,3\}$.
In the case $\chi(X)=1$ there may or may not be base points; in the case $\chi(X)=3$ there are always base points.
\end{thm}

Both results are sharp in the following sense: $|4I(K_X+\Delta)|$ is always base-point-free by \cite[Theorem~24]{LR14} and in general there is a local obstruction to the base-point-freeness of the bicanonical map and for $|3I(K_X+\Delta)|$ and $|4I(K_X+\Delta)|$  to be an embedding \cite[Example 46]{LR14}, which does not depend on the global volume of the surface.

\subsection*{Acknowledgements} We are grateful to Ben Antes, Fabrizio Catanese,  Stephen Coughlan, Marco Franciosi, Rita Pardini
for discussions about this problem in particular and algebraic surfaces in general. Wenfei Liu was supported by the NSFC (no.~12571046). Sönke Rollenske was supported by the DFG and thanks Xiamen University for the hospitality. We acknowledge the financial support for this collaboration provided by the Tianyuan Mathematical Center in Southeast China.

\subsection*{AI use disclosure} ChatGPT \cite{ChatGPT} puzzled together the proof from a rough set of notes accumulated over the years, suggesting in particular the positivity estimate in the proof of Lemma~\ref{lem: tricanonical curve} as an improvement to \cite[Lemma~40]{LR14} in the relevant special case, as well as the doubling construction to extend the proof to log surfaces. The paper was then jointly edited by ChatGPT, Claude \cite{Claude} and the authors until the latter were satisfied with the presentation.

\section{Preliminaries}
\subsection{Stable log surfaces}\label{sec: stable log surface}
We work over the complex numbers and refer to \cite{KollarSMMP, KollarModuli} for the general theory of slc pairs and their moduli. A \emph{demi-normal} surface is a reduced surface which satisfies Serre's condition $S_2$ and has at worst ordinary double points in codimension one. A \emph{stable log surface} $(X,\Delta)$ consists of a connected projective demi-normal surface $X$ and a reduced Weil divisor $\Delta$, none of whose components is contained in the non-normal locus of $X$, such that $(X,\Delta)$ has semi-log-canonical (slc) singularities and $K_X+\Delta$ is ample. If $\Delta=0$ we call $X$ a \emph{stable surface}. The \emph{global (Cartier) index} of $(X,\Delta)$ is the smallest positive integer $I$ such that
\[
 L=I(K_X+\Delta)
\]
is Cartier. Stable surfaces of index one are called Gorenstein. We write $\chi(X)=\chi(\mathcal O_X)$.

\subsection{Normalisation and gluing}\label{sec: gluing}
Let $\pi\colon\bar X\to X$ be the normalisation, let $D\subset X$ be the non-normal locus and $\bar D=\pi^{-1}(D)$ the conductor divisor, both with their reduced structure, and let $\bar\Delta$ be the strict transform of $\Delta$. Then $(\bar X,\bar D+\bar\Delta)$ is a log canonical pair, possibly disconnected, and
\[
 \pi^*(K_X+\Delta)=K_{\bar X}+\bar D+\bar\Delta.
\]
Since $X$ has ordinary double points at the generic points of $D$, the induced map of normalisations $\bar D^\nu\to D^\nu$ is the quotient by an involution $\tau$, and $\tau$ preserves the different of $(\bar X,\bar D+\bar\Delta)$ on $\bar D^\nu$:
\begin{equation*}
	\begin{tikzcd}
		\bar X \dar{\pi}\rar[hookleftarrow]{\bar\iota} & \bar D\dar{\pi} & \bar D^\nu \lar[swap]{\bar\nu}\dar{/\tau}
		\\
		X\rar[hookleftarrow]{\iota} &D &D^\nu\lar[swap]{\nu}
	\end{tikzcd}
\end{equation*}
By Koll\'ar's gluing theorem \cite[Theorem~5.13]{KollarSMMP}, $(X,\Delta)$ is uniquely determined by the triple $(\bar X,\bar D+\bar\Delta,\tau)$. Conversely, every log canonical pair with ample log canonical divisor, together with a generically fixed-point-free involution of the normalisation of a reduced part of its boundary that preserves the different, arises in this way.

Sections of the pluri-log-canonical bundles can be computed on the normalisation \cite[Proposition~5.8]{KollarSMMP}: for an integer $k$, a section $s$ of $\omega_{\bar X}(\bar D+\bar\Delta)^{[k]}$ over $\pi^{-1}(U)$ is the pullback of a section of $\omega_X(\Delta)^{[k]}$ over $U$ if and only if
\begin{equation}\label{eq: descent}
 \tau^*\Res_{\bar D^\nu}(s)=(-1)^k\Res_{\bar D^\nu}(s),
\end{equation}
where $\Res_{\bar D^\nu}$ is the residue map to the $k$-canonical forms on $\bar D^\nu$ with poles along the different. In particular, sections vanishing along $\bar D$ always descend; compare \cite[Corollary~20]{LR14}.

\subsection{slc centres and vanishing}\label{sec: vanishing}
The \emph{slc centres} of $(X,\Delta)$ are the images of the log canonical centres of $(\bar X,\bar D+\bar\Delta)$, and the \emph{slc strata} are the slc centres together with the irreducible components of $X$ \cite{Fujino14}. The one-dimensional slc centres are the components of $D$ and of $\Delta$. The zero-dimensional slc centres form a finite set, which contains the singular points of $D\cup\Delta$ \cite[Lemma~29]{LR14} and the normal points of $X\setminus\Delta$ at which $X$ is not klt.

Let $(X,\Delta)$ be a stable log surface of index $I$ and recall that $L = I(K_X+\Delta)$ is Cartier. Let $W\subset X$ be a union of slc strata with its reduced structure, and let $m\geq1$ with $mI\geq2$. Since $mL-(K_X+\Delta)=(mI-1)(K_X+\Delta)$ is ample, Fujino's vanishing theorems \cite[Theorems~1.7 and~1.11]{Fujino14} give
\begin{equation}\label{eq: Fujino vanishing}
 H^i(X,\mathcal O_X(mL))=0
 \quad\text{and}\quad
 H^i(X,\mathcal I_W\otimes\mathcal O_X(mL))=0
 \qquad\text{for } i>0.
\end{equation}
We will mostly use the following consequences:
\begin{enumerate}
 \item the restriction map $H^0(X,mL)\to H^0(W,mL|_W)$ is surjective;
 \item in particular, taking $W$ to be a point, $mL$ is globally generated at every zero-dimensional slc centre.
\end{enumerate}
Moreover, $W$ is seminormal \cite[Theorem~1.13]{Fujino14}.
%

\subsection{Very ampleness and curves}\label{sec: curves}
Let $H$ be a line bundle on a projective scheme $Y$ and $V\subset H^0(Y,H)$ a linear subspace. We say that $V$ (or the linear system $|V|$) \emph{separates} a subscheme $\xi\subset Y$ of length two if $V\to H^0(\xi,H|_\xi)$ is surjective. Thus $H$ is very ample if and only if $|H|$ separates every length-two subscheme $\xi$ of $Y$. This is usually split into separating points, when $\xi$ is two points, and separating tangent vectors, when $\xi$ is non-reduced supported at one point. The latter can also be expressed by saying that  $V$ generates first jets, that is, for every point $p$ with maximal ideal $\gothm_p$ the restriction map \[V \into H^0(Y, H) \to H^0(Y, \ko_Y/\gothm_p^2(H))\] is surjective.

It is often convenient to prove separation for a subscheme $\xi$ by restriction to a curve through $\xi$. To be precise, by a \emph{curve} we mean a projective Cohen--Macaulay scheme of pure dimension one, possibly reducible and nonreduced; for instance, every effective Cartier divisor on a stable surface is a curve. A \emph{subcurve} $B\subseteq C$ is a closed subscheme which is itself a curve, and $p_a(B)=1-\chi(\mathcal O_B)$, also when $B$ is disconnected. We use the following criterion \cite[Theorem~1.1]{CFHR}: a line bundle $H$ on a curve $C$ is very ample if
\begin{equation*}
 \deg(H|_B)\geq2p_a(B)+1
\end{equation*}
for every generically Gorenstein subcurve $B\subseteq C$, including $B=C$; it is globally generated if $\deg(H|_B)\geq2p_a(B)$ for all such $B$. If $C$ is reduced, every subcurve is reduced and hence generically Gorenstein.

\section{The non-Gorenstein case}
In this section we prove the following.
\begin{thm}\label{thm: high embedding}
Let $X$ be a stable surface of global index $I\geq2$. Then $mIK_X$ is very ample for every $m\geq5$.
\end{thm}
For $m\geq6$ this is \cite[Theorem~35$(iii)$]{LR14}, so it remains to show that $5IK_X$ is very ample. 

Let us set up some notation: we are considering the line bundle  $L=IK_X$ on a non-Gorenstein stable surface $X$, so  $I\geq 2$. 
For a component $\bar X_0$ of the normalisation put
\[A=(\pi^*K_X)|_{\bar X_0}=K_{\bar X_0}+\bar D_0,\] where $\bar D_0$ denotes the part of the conductor divisor on $\bar X_0$. Since $A$ is an integral Weil divisor
and $IA$ is ample Cartier,
\begin{equation*}
 IA^2=(IA)\cdot A\in\IZ_{>0},\qquad A^2\geq1/I.
\end{equation*}

We will construct sections of $mL$ on partial resolutions of the normalisation.
Let $f\colon S\to \bar X_0$ be a partial resolution, which is to be specified in different situations below, choose compatible
canonical divisors and write
$f^*A=K_S+\bar D_S+E$, where $\bar D_S$ is the strict transform of
$\bar D_0$ and $E$ is exceptional. For
$M_m=(mI-1)f^*A$, the Cartier divisor $mIf^*A$ gives
\begin{equation}\label{eq: rounded}
 K_S+\lceil M_m\rceil
 =mIf^*A-\bar D_S-\lfloor E\rfloor.
\end{equation}
The unrounded $\IQ$-divisor $M_m$ is nef and big. The global sections of the divisor $mIf^*A-\bar D_S-\lfloor E\rfloor$ push forward to global sections of 
$\mathcal O_{\bar X_0}(mIA-\bar D_0)$. They vanish along the conductor, so they descend to $X$ by \eqref{eq: descent}, extending by zero on the other components of $X$. Thus, if the map $S\rightarrow \bar X_0\rightarrow X$ is an isomorphism over a neighbourhood $U\subset X$, the global generation (resp.~very ampleness) of the divisor $ K_S+\lceil M_m\rceil$ over $U$ implies that of $mL$ over $U$.

\begin{lem}\label{lem: base locus}
For $I\geq 2$ the base locus of $|2L|$ is contained in the smooth locus of $X$, that is, $\Bs|2L|\subset X_{\mathrm{sm}}$.
\end{lem}
\begin{proof}
We will first treat the non-normal locus $D\subset X$, for which we can use Fujino's lifting results \eqref{eq: Fujino vanishing}. Then we handle the normal singular points by another approach. 

First we prove $L\cdot B\geq p_a(B)$ for every nonempty subcurve
$B\subset D$, which is reduced, because $D$ is reduced. Let $B_1,\ldots,B_a$ be the irreducible components of $B$, and let $g_j$ be the genus of the normalisation of $B_j$.
By \cite[Lemma~29]{LR14}, the curve $D$ has only multinodes, that is, it locally looks like the union of the coordinate axes in an affine space. Hence so does $B$, and a singular point $x$ of $B$ with $\mu_x$ branches contributes $\mu_x-1$ to the arithmetic genus. Let $N$ be the number of singular points of $B$, and let $r_j$ be the number of branches of $B_j$ passing through them. Then
\begin{equation}\label{eq: pa of B}
 p_a(B)=1+\sum_{j=1}^a(g_j-1)+\sum_{x\in\operatorname{Sing}B}(\mu_x-1)
 =1-N+\sum_{j=1}^a(g_j+r_j-1).
\end{equation}
On the other hand, \cite[Lemma~30]{LR14} and the ampleness of the Cartier divisor $L$ give
\[
 L\cdot B_j\geq I(2g_j-2+r_j)\geq 2(2g_j-2+r_j)\quad\text{and}\quad L\cdot B_j\geq1;
\]
in \cite{LR14} the branches are counted over all singular points of $D$, which only makes the first inequality stronger. From these two bounds we get, for every $j$,
\begin{equation}\label{eq: LBj}
 L\cdot B_j\geq g_j+r_j-1,\qquad\text{and}\quad L\cdot B_j\geq g_j\ \text{ if } r_j=0.
\end{equation}
Indeed, if $g_j\geq2$ then $2(2g_j-2+r_j)\geq g_j+r_j$; if $g_j=1$ then $2r_j\geq r_j$, and $L\cdot B_j\geq1$ covers $r_j=0$; if $g_j=0$ we need $L\cdot B_j\geq r_j-1$, which follows from $L\cdot B_j\geq1$ for $r_j\leq2$ and from $2(r_j-2)\geq r_j-1$ for $r_j\geq3$.

Now sum \eqref{eq: LBj} over $j$ and compare with \eqref{eq: pa of B}. If $B$ is singular, then $N\geq1$ and
\[
 L\cdot B\geq\sum_j(g_j+r_j-1)=p_a(B)+N-1\geq p_a(B).
\]
If $B$ is smooth, then all $r_j=0$, and $L\cdot B\geq\sum_j g_j\geq p_a(B)=1-a+\sum_{j=1}^a g_j$.

Thus, we have $\deg 2L|_B\geq 2p_a(B)$ for every nonempty reduced subcurve
$B\subset D$, and by the criterion in Section~\ref{sec: curves}, $2L|_D$ is globally generated.
By \eqref{eq: Fujino vanishing}, the sections of $2L|_D$ on $D$ lift to sections of $2L$, since $D$ is a union of slc centres. Hence $|2L|$ has no base points along $D$.

It remains to consider a singular point of $X$ outside $D$, that is, the image $\pi(p)$ of a singular point $p\in \bar X_0\setminus\bar D_0$. Let $f\colon S\rightarrow \bar X_0$ resolve all other singularities and
retain $p$, so that $S\rightarrow \bar X_0\rightarrow X$ is an isomorphism over a neighbourhood of $\pi(p)$. The divisor $K_S+\lceil M_m\rceil$ of \eqref{eq: rounded} is Cartier on $S$:
away from $f^{-1}(p)$ the surface is smooth, and near $f^{-1}(p)$ it equals $mIf^*A$.
The fractional part $\lceil M_m\rceil-M_m$ zero in a neighbourhood of $p$.
For every integral curve $B$ through $p$,
\begin{equation}\label{eq: M2}
 M_2^2\geq(2I-1)^2/I\geq9/2,
 \qquad M_2\cdot B\geq(2I-1)/I\geq3/2.
\end{equation}
Since the boundary $\lceil M_2\rceil-M_2$ does not pass through $p$, Kawachi's freeness theorem \cite[Theorem~1]{Kawachi98} shows that $p$ is not a base point of $|K_S+\lceil M_2\rceil|$ as soon as $M_2^2>\delta_p$ and $M_2\cdot B\geq\delta_p/2$ for every curve $B$ through $p$. Here $\delta_p\leq 2$ if $p$ is klt \cite[Proposition~1]{Kawachi98} and $\delta_p=0$ otherwise. The inequalities \eqref{eq: M2}
therefore give global generation of $2IA-\bar D_0$ at $p$, which descends to global generation of $2L$ of $X$ at $\pi(p)$.
\end{proof}

\begin{lem}\label{lem: smooth separation}
	For $I\geq 2$ the linear system $|5L|$ separates every subscheme of length
two supported in the smooth locus of $X$.
\end{lem}
\begin{proof}
Let $\xi$ be supported at the smooth points $p,q$ of $X$, where $p=q$ is allowed.

Assume first that $p$ and $q$ lie on the same component $X_0$ of $X$, dominated by the normalisation component $\bar X_0$, and identify $p$, $q$ with their preimages in $\bar X_0$. Take a resolution $f\colon S\rightarrow \bar X_0$ which is an isomorphism over neighbourhoods of $p$ and $q$. Then
\[
 M_5^2\geq(5I-1)^2/I\geq81/2>8,
 \qquad M_5\cdot B\geq5-1/I\geq9/2>4
\]
for every integral curve $B$ through $f^{-1}(p)$ or $f^{-1}(q)$, and the fractional
part $\lceil M_5\rceil-M_5$ is zero in a neighbourhood of these points. The local separation theorems of
\cite[Proposition~4(3) and Proposition~5, Case~3]{Masek99} apply with
all parameters $\beta_2=\beta_1=2$, which correspond to the conditions
$M_5^2>8$ and $M_5\cdot B\geq4$. They separate both two distinct
points and an infinitesimal tangent vector at a point.

The sections of $K_S+\lceil M_5\rceil$ descend to $X$ as above, so $|5L|$ separates $\xi$.

If $p$ and $q$ lie on distinct components of $X$, then the same argument on the component of $q$ gives a section of $5L$ on that component which is nonzero at $q$; extended by zero to the other components it gives $s\in H^0(X, 5L)$ with $s(p)=0$ and $s(q)\neq 0$. Symmetrically there is $t\in H^0(X,5L)$ with $t(p)\neq 0$ and $t(q)= 0$, so $|5L|$ separates $p$ and $q$.
\end{proof}

\begin{lem}\label{lem: high curve}
Assume $I\geq 2$ and let $C\in|3L|$ be reduced, containing no component of the non-normal locus $D$. Then $5L|_C$ is very ample, and all its sections lift to $X$.
\end{lem}
\begin{proof}
The curve $C$ is Cartier on $X$ and hence is Cohen--Macaulay. For a nonempty reduced
subcurve $B\subset C$, put $d=L\cdot B\in\IZ_{\geq1}$.
By \cite[Lemma~40]{LR14},
\[
 2p_a(B)-2\leq(3+1/I)d.
\]
If $d=1$, the even integer on the left is at most $7/2$, hence at
most two; thus $p_a(B)\leq2$ and $5d\geq2p_a(B)+1$.
For $d\geq2$, the surplus is
\[
 5d-(2p_a(B)-2)\geq(2-1/I)d\geq3.
\]
Applying \cite[Theorem~1.1]{CFHR} gives the very ampleness of $5L|_C$. 

The long exact sequence associated to 
$0\to\mathcal O_X(2L)\to\mathcal O_X(5L)\to\mathcal O_C(5L)\to0$
and the vanishing \eqref{eq: Fujino vanishing} prove the lifting of the sections of $5L|_C$ to $X$.
\end{proof}

\begin{proof}[Proof of Theorem~\ref{thm: high embedding}]
For $m\geq6$ the assertion is \cite[Theorem~35$(iii)$]{LR14}, so it suffices to show that $|5L|$ separates every length-two subscheme $\xi\subset X$. The linear system $|3L|$ is base-point-free by
\cite[Theorem~24]{LR14}. If it does not separate $\xi$,
\cite[Lemma~38$(i)$]{LR14} supplies a reduced Cartier curve
$\xi\subset C\in|3L|$ containing no component of the non-normal locus $D$. By
Lemma~\ref{lem: high curve}, $5L|_C$ separates $\xi$ and its sections lift to $X$. Hence $|5L|$ embeds $\xi$.

Suppose now that $|3L|$ separates $\xi$. If the support of $\xi$ avoids $\Bs|2L|$, then choose
a section $s$ of $2L$, nonzero at every point of the support of  $\xi$. Then the sublinear system $(s=0)+|3L|\subset |5L|$ and hence also $|5L|$ separates $\xi$.

If $\xi=\{p,q\}$ with $p\in\Bs|2L|$ and $q\notin\Bs|2L|$, choose $s_2\in H^0(X,2L)$ and $s_3\in H^0(X,3L)$ not vanishing at $q$; then $s=s_2s_3$ vanishes at $p$ but not at $q$. Since $|5L|$ is base-point-free \cite[Theorem~24]{LR14}, there is $u\in H^0(X,5L)$ with $u(p)\neq0$, and a suitable linear combination of $u$ and $s$ vanishes at $q$ but not at $p$. Thus $|5L|$ separates $p$ and $q$.

In all remaining cases the support of $\xi$ lies in
$\Bs|2L|$, which is contained in $X_{\mathrm{sm}}$ by Lemma~\ref{lem: base locus}, and $|5L|$
separates $\xi$ by Lemma~\ref{lem: smooth separation}. 
\end{proof}

\section{The Gorenstein case}
\subsection{Exceptional cases  for the tricanonical map}

The following result is classical for surfaces of general type \cite[VII.7]{BHPV}. A problem we have to face in the Gorenstein stable case is that, even when the global volume is large, we can have irreducible components of volume one. 
\begin{thm}\label{thm: gorenstein tricanonical}
	Let $X$ be a Gorenstein stable surface.
	\begin{enumerate}
		\item If $K_X^2 = 1$ and $\chi(X) = 1$, then $|3K_X|$ can have base points.
		\item If $K_X^2 = 1$ and $\chi(X) = 3$, then $|3K_X|$ has a base point. 
		\item In all other cases $|3K_X|$ is base-point-free. 
	\end{enumerate}
\end{thm}
Since the surfaces of volume one play a special role, we list them in Table~\ref{tab: volume one}. Especially the ones with $\chi(X) = 3$, sometimes called I-surfaces or $(1,2)$-surfaces, have attracted a lot of attention over the last years \cite{FPRR22, CFPRR, GPSZ22, RollenskeTorres, CFPR-Hodge, CFPR23, CFPR-2Gor}. They play a role in characterising threefolds of general type on the Noether line \cite{CP23, CHPZ24}. 

The classification of numerical Godeaux surfaces, $K_X^2 = \chi(X) = 1$,  is a longstanding problem and despite recent progress, their classification is as yet incomplete \cite{MR4771342, MR5065701}. The Gorenstein stable ones with worse than canonical singularities have been classified in \cite{FPRGodeaux}. 
\begin{table}
\begin{tabular}{ccccc}
\toprule
$\chi(X)$ & $p_g(X)$ & $q(X)$ & description & reference\\
\midrule
3 & 2 & 0 & $X_{10}\subset\IP(1,1,2,5)$ & \cite{FPR17}\\
2 & 1 & 0 & $X_{6,6}\subset\IP(1,2,2,3,3)$ & \cite{FPR17}\\
1 & 0 & 0 & stable Godeaux surfaces & \cite{FPRGodeaux}\\
0 & 0 & 1 & $X_{0.1}$, $X_{0.2}$ & see Section~\ref{sec: computations}\\
\bottomrule
\end{tabular}
\caption{Gorenstein stable surfaces with $K_X^2=1$.}\label{tab: volume one}
\end{table}

\begin{proof}
The first two  items follow easily from classical results and the  descriptions in Table~\ref{tab: volume one}.

For a smooth numerical Godeaux surface, Miyaoka's theorem describes the
base points of $|3K_X|$ in terms of its torsion group
\cite[Theorems~2 and~2']{MR0409481}.  Both possibilities in case~$(i)$
occur already among these surfaces.

As observed in  \cite[Proposition~3.6$(iii)$]{FPR17} or classically in  \cite[VII.7]{BHPV} the base point of the canonical pencil of $X_{10} \subset \IP(1,1,2,5)$ is also a base point of the tricanonical map, which proves $(ii)$.

To prove $(iii)$ we need to show that $|3K_X|$ has no base points in all other cases.
We first exclude the remaining cases of volume one: for $K_X^2=1$ and $\chi(X)=2$, generation of $3K_X$ is \cite[Proposition~3.7$(ii)$]{FPR17}. The two surfaces with $\chi(X)=0$ need to be treated by a direct computation, which we relegate to Section~\ref{sec: computations}, 
Proposition~\ref{prop: chi zero surfaces} in order not to break the flow of the proof. 

Assume henceforth that $K_X^2>1$. 
The global generation of $3K_X$ along the non-klt locus of $X$ is verbatim Steps~2--3 of \cite[proof of Theorem~5.1]{Fuj17} with $D=3K_X$ and $A=D-K_X=2K_X$. 
%

Hence the base locus of $3K_X$ is contained in the normal locus, and we treat this component by component.

On a normalisation component $\bar X_0\subset \bar X$ with conductor divisor $\bar D_0$ write
\[
 A=\pi^*K_X|_{\bar X_0}=K_{\bar X_0}+\bar D_0.
\]
It is ample and Cartier, so $A^2$ is a positive integer.  If $A^2>1$,
\cite[Proposition~26]{LR14}, with $m=3$, produces sections of $3A$
vanishing along $\bar D_0$ and generating it outside $\bar D_0$.
Extended by zero to the other normalisation components, these sections descend to $X$. Thus $3K_X$ is generated by global sections also on the normal locus of any component of volume at least two.

It remains to consider $A^2=1$.  Let $X_0$ be the reduced irreducible
component of $X$ dominated by $\bar X_0$, so that $\pi_0=\pi|_{\bar X_0}\colon \bar X_0\to X_0$ is its normalisation.
Since $K_X^2>1$, $X_0$ is a proper component of $X$.  Connectedness and the $S_2$
property imply that $X_0$ meets the other components along a curve.
Therefore
we can decompose the reduced conductor on $\bar X_0$ as
\[
 \bar D_0=\bar\Gamma+\bar B,\qquad \bar B\ne0,
\]
where $\bar\Gamma$ denotes the inverse image of the internal non-normal locus on $X_0$,
and $\bar B$ denotes the inverse image of the intersection of $X_0$ with other components of $X$.

The component $X_0$ is an slc stratum, hence is seminormal by
\cite[Theorem~1.13]{Fujino14}.  Let $\mathfrak c_0$ be the conductor
ideal of $\pi_0$, considered on $\bar X_0$.
The conductor of a finite normalisation of a seminormal scheme is radical \cite[Corollary~1.5(1)]{LV81}.  

The Gorenstein gluing description \cite[Theorem~3.2 and its
addendum]{FPR15} shows that the support of this conductor is contained
in the union of $\bar\Gamma$ and the nodes of $\bar D_0$.  Indeed, outside $\bar D_0$ the normalisation $\pi_0$
is an isomorphism, and at a smooth point of $\bar B$, $X_0$ is already normal.  Singularities of $\bar X_0$ along
$\bar D_0$ occur only at its nodes.  Consequently
\begin{equation}\label{eq: descend to X0}
 \mathcal I_{\bar\Gamma\cup\operatorname{Sing}(\bar D_0)}\subseteq\mathfrak c_0,
 \qquad
 (\pi_0)_*(\mathfrak c_0\otimes \pi_0^*\mathcal O_{X_0}(3K_X))
 \subseteq\mathcal O_{X_0}(3K_X).
\end{equation}
Thus sections of $3A$ vanishing on $\bar\Gamma$ and at the nodes of $\bar D_0$ descend to
$X_0$.  Finally,
\[
 H^0(X,3K_X)\longrightarrow H^0(X_0,3K_X|_{X_0})
\]
is surjective by \eqref{eq: Fujino vanishing}, since $X_0$ is an slc stratum. Thus it suffices to find, for every point $p\in\bar X_0\setminus\bar D_0$, a section of $3A$ which vanishes on $\bar\Gamma$ and at the nodes of $\bar D_0$, but not at $p$.

The possible pairs $(\bar X_0, \bar D_0)$ have been classified in \cite[Theorem~1.1 and List~2.2]{FPR15} and we treat them
case by case.

\begin{prooflist}
\item[The case $(P)$]
Here $\bar X_0=\IP^2$, $A=\mathcal O_{\IP^2}(1)$, and $\bar D_0$ is a
nodal quartic.  Write its equation as $\gamma\beta =0$ such that $\bar\Gamma = (\gamma=0)$ and $\bar B=(\beta=0)$; take $\gamma=1$ if $\bar\Gamma=0$.
The cubics
\[
 s_i=\gamma\frac{\partial\beta}{\partial x_i},\qquad i=0,1,2,
\]
vanish on $\bar\Gamma$ and at every node of $\bar D_0$:
At a node of $\bar B$ the derivatives $\frac{\partial\beta}{\partial x_i}$ vanish; at a node on
$\bar\Gamma$ the factor $\gamma$ vanishes. Thus the $s_i$ descend to sections of $\mathcal{O}_{X_0}(3K_X)$ on $X_0$ by \eqref{eq: descend to X0}.  On the complement of
$\bar D_0$, Euler's identity
\[
 \sum_{i=0}^2x_is_i=(\deg\beta)\gamma\beta
\]
is nonzero, so these sections generate $3A$ there, and they descend to $X_0$.

We claim that $\bar \Gamma=0$ in the other three cases $(dP)$, $(E_-)$, and $(E_+)$. In fact, if $\bar\Gamma\ne0$, the map from $\bar \Gamma$ to its image $\Gamma_X$ on $X$ has generic
degree two.  The projection formula gives
\[
 A\cdot \bar\Gamma=2K_X\cdot \Gamma_X\in2\mathbb Z_{>0},
 \qquad A\bar B\ge1.
\]
This contradicts the equality $A\cdot(\bar\Gamma+\bar B)=2$ given by the classification of \cite{FPR15}.

\item[The cases $(dP)$ and $(E_-)$]
Here $\bar D_0=\bar B$, $A=-K_{\bar X_0}$ and $\bar B\in|2A|$. In both cases $\bar X_0$ is a
Gorenstein del Pezzo surface of degree one, with canonical singularities in case $(dP)$ and with
an elliptic singularity in case $(E_-)$. Its anticanonical ring exhibits it as a sextic \cite[Theorem~4.4$(iv)$]{HW81}
\begin{equation*}
 \bar X_0=\{z^2=f_6(x_0,x_1,y)\}\subset\IP(1,1,2,3),
 \qquad A=\mathcal O_{\bar X_0}(1),
\end{equation*}
where $x_0,x_1,y,z$ have degrees $1,1,2,3$. The coefficient of $y^3$ in $f_6$ is nonzero, since the
point $(0:0:1:0)$, a non-Gorenstein singular point of $\IP(1,1,2,3)$, does not lie on $\bar X_0$ \cite[Proposition~4.6]{HW81}.

The curve $\bar B$ is the zero locus of a quadric
$s\in H^0(\bar X_0, 2A) = \langle x_0^2,x_0x_1,x_1^2,y\rangle$,
and we look for cubics vanishing at the nodes of $\bar B$ without common zero outside $\bar B$.
There are two possibilities.

If $s$ involves $y$, we may take $s$ as a new coordinate in place of $y$, so that $\bar B=\{y=0\}$.
Then $\bar B$ is the double cover $z^2=f_6(x_0,x_1,0)$ of $\IP^1$, and its nodes lie over the
multiple roots of the sextic $f_6(x_0,x_1,0)$; in particular, $z=0$ at every node.
Hence the cubics $x_0y$, $x_1y$ and $z$ vanish at all nodes of $\bar B$. A common zero outside
$\bar B$ would have $y\neq0$ and $x_0=x_1=z=0$, so it would be the point $(0:0:1:0)$,
which does not lie on $\bar X_0$.

If $s$ does not involve $y$, it is a quadric in $x_0,x_1$, and it is not a square because
$\bar B$ is reduced. After a linear change of $x_0,x_1$ we have $s=x_0x_1$, that is,
$\bar B=\{x_0=0\}\cup\{x_1=0\}$. The single cubic $x_0^2x_1$ vanishes on all of $\bar B$, in
particular at its nodes, and it vanishes nowhere outside $\bar B$.

\item[The case $(E_+)$]
Here $\bar X_0=S^2E$ is smooth, and, with the section and fibre classes
$C_0,F$ of \cite[List~2.2]{FPR15},
\[
 C_0^2=C_0\cdot F=1,\quad F^2=0,\quad
 A\equiv C_0,\quad K_{\bar X_0}\equiv-2C_0+F,\quad
 \bar B\equiv3C_0-F.
\]
The divisor $-K_{\bar X_0}$ is nef: its pullback under the finite quotient
$E\times E\to S^2E$ is the diagonal, a nef divisor on the abelian
surface.  Since $(-K_{\bar X_0})\cdot \bar B=1$ and $(-K_{\bar X_0})\cdot F=2$, the reduced
divisor $\bar B$ has no fibre component, so every component of $\bar B$ dominates $E$ and its
normalisation has genus at least one. By adjunction, $\bar B$ is a nodal curve of arithmetic genus
$1+\frac12 A\cdot\bar B=2$. If it has $r$ components and $\delta$ nodes, then
\[
 2=\sum_{j=1}^r g(\bar B_j^\nu)+\delta-r+1\ge\delta+1,
\]
so $\bar B$ has at most one node. By the gluing description, points of $\bar X_0$ can only be
identified with each other if they lie over nodes of $\bar D_0$, so with at most one node
nothing is glued: $\pi_0$ is bijective, hence an isomorphism since $X_0$ is seminormal.
It therefore suffices to show that $3A$ is globally generated.

Put $N=3A-K_{\bar X_0}$.  It is nef and satisfies
\[
 N^2=9A^2-6A\cdot K_{\bar X_0}+K_{\bar X_0}^2=15
 \text{ and } N\cdot C=3A\cdot C+(-K_{\bar X_0})\cdot C\ge3
\]
for every integral curve $C$, because $A$ is ample and Cartier.
Reider's theorem \cite[Theorem~1]{Reider88} implies that
$K_{\bar X_0}+N=3A$ is very ample, in particular globally generated.
\end{prooflist}
The proof of $(iii)$ is now complete.
\end{proof}

\subsection{Very ampleness of \texorpdfstring{$mK_X$}{mK\_X} in the Gorenstein case}
We now prove the extension of the classical bound by Bombieri \cite{bombieri73} to the Gorenstein stable case. 
\begin{thm}\label{thm: very ample}
Let $X$ be a Gorenstein stable surface. Then $mK_X$ is very ample for every $m\geq5$.
\end{thm}

Separating points will be relatively easy, the main problem is the separation of tangent vectors. 
At the klt points of $X$, which in the Gorenstein case means smooth points and isolated canonical singularities, we have powerful Reider-type results at our disposal and for subschemes contained in the non-klt locus we have the restriction results by Fujino \eqref{eq: Fujino vanishing}. 

The main difficulties are presented by subschemes supported on but not contained in the non-klt locus. Those on semi-smooth points can be treated by hand, Proposition~\ref{prop: local separation}, for the other ones we lift sections from suitable bi- or tri-canonical curves.

Throughout this subsection, $X$ is a Gorenstein stable surface with normalisation $\pi\colon\bar X\to X$, non-normal locus $D$ and conductor $\bar D$, and $m\geq5$ is an integer. Let $\Sigma\subset X$ be the finite set of zero-dimensional slc centres, that is, the singular points of $D$ and the normal points at which $X$ is not klt (Section~\ref{sec: vanishing}), and let $\Nklt(X)=D\cup\Sigma$ be the non-klt locus, with its reduced structure.

Recall that, by the classification of Gorenstein slc surface singularities \cite[Remark~2.9]{LR16}, a klt point of $X$ is a smooth point or a rational double point, and a point of $D\setminus\Sigma$ is a double crossing point or a pinch point.

The first step is the following. 
\begin{prop}\label{prop: local separation}
Let $X$ be a Gorenstein stable surface and $m\geq5$.
Then $|mK_X|$ separates every pair of distinct points on $X$ and every length-two
subscheme whose support lies in $X\setminus\Sigma$.
\end{prop}
We will collect several auxiliary results before we embark on the proof. 

\begin{lem}\label{lem: resolution}
Let $g\colon S\to\bar X$ be the minimal resolution of the singularities of $\bar X$ lying over $\Sigma$,
\[
\begin{tikzcd}
  S\rar{g}\drar[swap]{f} & \bar X\dar{\pi}\\
  & X
 \end{tikzcd},
 \]
 and put $A=f^*K_X$ and $M=(m-1)A$.
 In this situation:
\begin{enumerate}
 \item $S$ has at most rational double points, and $f$ is finite over $X\setminus\Sigma$ and an isomorphism over $X\setminus \Nklt(X)$. Near the preimages of a point of $D\setminus\Sigma$, the surface $S$ is smooth and $K_S+M=f^*(mK_X)-\bar D$.
 \item $A\cdot B\geq1$ for every curve $B\subset S$ meeting $f^{-1}(X\setminus\Sigma)$.
 \item Every section of $K_S+M$ descends to a section of $mK_X$ vanishing on the non-klt locus of $X$. 
 \end{enumerate}
\end{lem}
\begin{proof}
The first item follows from the description of the singularities of $X\setminus\Sigma$ above, together with $K_{\bar X}=\pi^*K_X-\bar D$. For the second, $f(B)$ is a curve because $f$ is finite over $X\setminus\Sigma$, so $A\cdot B$ is a positive multiple of the integer $K_X\cdot f(B)>0$.

For the third, recall that $g_*\mathcal O_S(K_S)\subseteq\mathcal O_{\bar X}(K_{\bar X})$ because $\bar X$ is normal. Since $K_{\bar X}=\pi^*K_X-\bar D$, this gives $g_*\mathcal O_S(K_S+M)\subseteq\mathcal O_{\bar X}\bigl(\pi^*(mK_X)-\bar D\bigr)$, so every section of $K_S+M$ is a section of $\pi^*(mK_X)$ vanishing along $\bar D$. Such sections descend to $X$ by \eqref{eq: descent}, and the descended section vanishes along $D$. Now let $p\in\Sigma\setminus D$. Then $p$ is a Gorenstein singularity which is not a rational double point, so near $g^{-1}(p)$ we have $K_S=g^*K_{\bar X}-Z_K$, where the canonical cycle $Z_K$ is a nonzero effective divisor supported on $g^{-1}(p)$ \cite[4.19--4.20]{ReidChapters}. Hence $g_*\mathcal O_S(K_S)\subseteq\mathfrak m_p\,\mathcal O_{\bar X}(K_{\bar X})$ near $p$, and the descended sections vanish at $p$ as well.
\end{proof}

\begin{lem}\label{lem: S separation}
The linear system $|K_S+M|$ separates any two distinct points of $f^{-1}(X\setminus\Sigma)$ and generates first jets at every point of $f^{-1}(X\setminus\Sigma)$.
\end{lem}
\begin{proof}
We work on a connected component of the normal surface $S$, retaining the notation $S$, $A$ and $M$ for the component and the restricted divisors. We apply \cite[Theorems~1, 3 and~4]{KM98} to the nef divisor $M$. Since $A^2$ is a positive integer, we have $M^2=(m-1)^2A^2\geq16$, and $M\cdot B=(m-1)A\cdot B\geq4$ for every curve $B$ through a point of $f^{-1}(X\setminus\Sigma)$ by Lemma~\ref{lem: resolution}$(ii)$. By Lemma~\ref{lem: resolution}$(i)$, these points are smooth points or rational double points of $S$, where the invariant $\delta_p$ of \cite{KM98} is at most $4$ \cite[Theorem~2]{KM98}. Hence any two such points on the same component are separated by \cite[Theorem~3]{KM98}; points on different components are separated by choosing sections independently on the components.

For first jets at $p$ we use \cite[Theorem~4]{KM98}. If $p$ is a rational double point, then $\delta_p(1)=8$ \cite[p.~247]{KM98}, and the bounds above suffice. If $p$ is a smooth point, then $\delta_p(1)=9$ \cite[p.~247]{KM98}. If first-jet generation failed at $p$, \cite[Theorem~4]{KM98} would give an effective divisor $B$ through $p$ with
\[
 0<M\cdot B<\frac92,\qquad
 M\cdot B-\frac94\leq B^2
 \leq\frac{(M\cdot B)^2}{M^2}=\frac{(A\cdot B)^2}{A^2}.
\]
Since $A\cdot B$ is a positive integer and $M\cdot B=(m-1)A\cdot B\geq4$, this forces $m=5$ and $A\cdot B=1$. The inequalities then give $7/4\leq B^2\leq1$, a contradiction.
\end{proof}

\begin{proof}[Proof of Proposition~\ref{prop: local separation}]
We use two sources of sections of $mK_X$.  First, $mK_X|_D$ is very ample by \cite[Proposition~31]{LR14}, and thus also on the full non-klt locus, because we only add a finite number of points in $\Sigma$. By \eqref{eq: Fujino vanishing} these sections lift to global sections on $X$.  In particular, $mK_X$ is globally generated along $\Nklt(X)$. Second, by Lemma~\ref{lem: resolution}$(iii)$ the sections of $K_S+M$ descend to sections of $mK_X$ vanishing on $\Nklt(X)$.

Let $\xi\subset X$ be a length-two subscheme which is either reduced or supported at a point of $X\setminus\Sigma$. Since $\Nklt(X)\setminus\Sigma=D\setminus\Sigma$, exactly one of the following cases occurs.
\begin{prooflist}
\item[$\xi\subset \Nklt(X)$] Then $\xi$ is separated, because $mK_X|_{\Nklt(X)}$ is very ample and sections lift by \eqref{eq: Fujino vanishing}.
\item[$\xi\cap \Nklt(X)=\emptyset$] Since $f$ is an isomorphism over $X\setminus \Nklt(X)$, the subscheme $\xi$ is the isomorphic image of a length-two subscheme of $f^{-1}(X\setminus \Nklt(X))$. The latter is separated by $|K_S+M|$, so $\xi$ is separated by the descended sections.
\item[$\xi=\{p,q\}$ with $p\notin \Nklt(X)$ and $q\in \Nklt(X)$] There is a descended section $s$ with $s(p)\neq0$, and $s(q)=0$ because $s$ vanishes on $\Nklt(X)$. Since $mK_X$ is globally generated along $\Nklt(X)$, there is a section $\sigma$ with $\sigma(q)\neq0$. Then $s$ and $\sigma-\frac{\sigma(p)}{s(p)}s$ separate $p$ and $q$.
\item[$\xi$ is supported at a point $p\in D\setminus\Sigma$ and $\xi\not\subset D$] We show that $|mK_X|$ generates all first jets at $p$, that is, the surjectivity of the middle row in the diagram
\begin{equation}\label{eq: jet diagram}
 \begin{tikzcd}[column sep=1.2em]
  && 0 \dar \\
  H^0(S, K_S+M) \rar[hookrightarrow]\arrow[bend left = 20] {rr}{\delta}  & H^0(X, \ki_D(mK_X))  \dar \rar & H^0(X,  \ki_D/(\ki_D\cap \gothm_p^2)(mK_X))    \dar\\
  & H^0(X, mK_X) \rar\dar[->>] & H^0(X, \ko_X/\gothm_p^2(mK_X))\dar \\
  & H^0(D, mK_X|_D) \rar[->>]&  H^0(D, \ko_X/( \ki_D + \gothm_p^2) (mK_X)) \dar
  \\
  &&0.
 \end{tikzcd}
\end{equation}
The surjectivity in the last row is \cite[Proposition~31]{LR14} and all the sections on $D$ lift by \eqref{eq: Fujino vanishing}. 
The map $\delta$ comes from Lemma \ref{lem: resolution} and we show its surjectivity using Lemma \ref{lem: S separation} now.

Let $\eta\subset S$ be the following length-two subscheme over $p$: if $p$ is a double crossing point, then $\eta=\{p_1,p_2\}$ consists of the two preimages of $p$; if $p$ is a pinch point, then $\eta$ is the tangent vector to $\bar D$ at the unique preimage $\tilde p$ of $p$. Near $\eta$ we have $K_S+M=\pi^*(mK_X)-\bar D$ by Lemma~\ref{lem: resolution}$(i)$, and the conductor $\ki_D$ is also the ideal of $\bar D$ in $\pi_*\ko_{\bar X}$. So, if $\ki_D\cap\gothm_p^2=\ki_D\cdot\ki_\eta$, then $\delta$ is the restriction map
\[
 H^0(S,K_S+M)\longrightarrow H^0\bigl(\eta,\ko_\eta(K_S+M)\bigr),
\]
which is surjective by Lemma~\ref{lem: S separation}. The equality can be read off from the local models, where $(x,y,z)=(uv,v,u^2)$ for the pinch point:
\begin{center}
\begin{tabular}{lll}
\toprule
 & double crossing point & pinch point\\
\midrule
 $\hat{\ko}_{X,p}$ & $\IC[[x,y,z]]/(xy)$ & $\IC[[x,y,z]]/(x^2-y^2z)$\\
 normalisation & $\IC[[x,z]]\times\IC[[y,z]]$ & $\IC[[u,v]]$\\
 $\ki_D$ & $(x)\times(y)$ & $(v)$\\
 $\ki_\eta$ & $(x,z)\times(y,z)$ & $(u^2,v)$\\
 $\ki_D\cap\gothm_p^2$ & $(x^2,xz)\times(y^2,yz)$ & $(u^2v,v^2)$\\
\bottomrule
\end{tabular}
\end{center}

Since the right-hand column is exact, a diagram chase in \eqref{eq: jet diagram} now shows that the middle row is surjective.
\qedhere
\end{prooflist}
\end{proof}

\begin{lem}\label{lem: curve trick}
Let $X$ be a Gorenstein stable surface, $m\geq5$, $p\in\Sigma$, $\xi\subset X$ a length-two subscheme supported at $p$, and $k\in\{2,3\}$. Assume that, if $|kK_X|$ does not separate $\xi$, there is a curve $C\in|kK_X|$ containing $\xi$ such that $mK_X|_C$ is very ample. Then $|mK_X|$ separates $\xi$.
\end{lem}
\begin{proof}
Since $p$ is a zero-dimensional slc centre and $m-k\geq2$, there is a section $\sigma$ of $(m-k)K_X$ with $\sigma(p)\neq0$ by \eqref{eq: Fujino vanishing}. If $|kK_X|$ separates $\xi$, then so does $\sigma\cdot|kK_X|\subset|mK_X|$.
Otherwise let $C$ be as in the assumption. The cohomology sequence of
\[
 0\longrightarrow\mathcal O_X((m-k)K_X)
 \longrightarrow\mathcal O_X(mK_X)
 \longrightarrow\mathcal O_C(mK_X)\longrightarrow0
\]
and $H^1(X,(m-k)K_X)=0$ show that every section of $mK_X|_C$ lifts to $X$, so $|mK_X|$ separates $\xi\subset C$.
\end{proof}

It remains to produce the curves. We do this for the two values of $k$.

\begin{lem}\label{lem: tricanonical curve}
Let $X$ be a Gorenstein stable surface and let
$C\in|3K_X|$ be a reduced curve.  Then $mK_X|_C$ is very
ample for every $m\geq 5$.
\end{lem}
\begin{proof}
By adjunction $\omega_C\cong\mathcal O_C(4K_X)$, which is globally generated by \cite[Theorem~24]{LR14}, and $C$ is connected because $H^1(X,-3K_X)\cong H^1(X,4K_X)^*=0$ by \eqref{eq: Fujino vanishing}.

Let $B\subsetneq C$ be a proper nonempty subcurve, $B'$ the complementary subcurve, and $t$ the length of $T=B\cap B'$. Finite duality for $B\hookrightarrow C$ gives
\begin{equation}\label{eq: subcurve duality}
 \omega_B=\mathcal H\!om_C(\mathcal O_B,\omega_C)\simeq\mathcal I_{T,B}\otimes\omega_C|_B,
 \quad\text{hence}\quad
 4K_X\cdot B=2p_a(B)-2+t.
\end{equation}
We claim $t\geq2$. Since $C$ is connected, every connected component of $B$ meets $B'$, so $t\geq1$, and if $t=1$ then $B$ is connected and $T$ is a reduced point $p$. In that case $h^1(B,\omega_B)=1$, while $H^1(B,\omega_C|_B)^*=H^0(B,\mathcal I_{p,B})=0$. So in the sequence $0\to\omega_B\to\omega_C|_B\to\omega_C|_B\otimes k(p)\to0$ the connecting map is surjective, and every section of $\omega_C|_B$ vanishes at $p$. This contradicts the global generation of $\omega_C$.

Hence $mK_X\cdot B-(2p_a(B)-2)=(m-4)K_X\cdot B+t\geq3$, and for $B=C$ adjunction gives $mK_X\cdot C-(2p_a(C)-2)=3(m-4)K_X^2\geq3$. By the criterion in Section~\ref{sec: curves}, $mK_X|_C$ is very ample.
\end{proof}

\begin{lem}\label{lem: bicanonical curve}
Let $X$ be a Gorenstein stable surface with $K_X^2=1$ and
$p_g(X)=0$.  Suppose that $C\in |2K_X|$ is not supported in $\Bs|2K_X|$. Then $mK_X|_C$ is very ample for every $m\geq5$.
\end{lem}
\begin{proof}
By the criterion in Section~\ref{sec: curves}, to show the very ampleness of $mK_X|_C$, it suffices to check that $mK_X\cdot B\geq 2p_a(B)+1$ for any generically Gorenstein subcurve $B$ of $C$.

Since $X$ is Gorenstein, $(K_X^2|_{X_i})$ is a positive integer for any irreducible component $X_i$ of $X$. Thus $K_X^2=1$ implies that $X$ is irreducible.
Adjunction gives
\[
 \omega_C\simeq\mathcal O_C(3K_X),\qquad p_a(C)=4,
 \qquad K_X\cdot C=2.
\]
Consequently $mK_X\cdot C=2m\geq10\geq2p_a(C)+1=9$.

If $C$ is reduced and irreducible, there are no proper nonempty
subcurves.  If it is reduced and reducible, it has exactly two
irreducible components $B,B'$, with $K_X\cdot B=K_X\cdot B'=1$. The vanishing $H^1(X,-2K_X)\simeq H^1(X,3K_X)^*=0$ shows that $C$ is connected, as in the proof of Lemma~\ref{lem: tricanonical curve}. The finite
duality argument in \eqref{eq: subcurve duality} gives, with
$t=\operatorname{length}(B\cap B')\geq 1$,
\[
 2p_a(B)-2=3K_X\cdot B-t\leq2.
\]
The left-hand side is even, so $p_a(B)\leq2$ and
$mK_X\cdot B=m\geq2p_a(B)+1$; the same holds for $B'$.

If $C$ is nonreduced, then $C=2B$, where
$B=C_{\mathrm{red}}$ is irreducible and $K_X\cdot B=1$. Its only nonempty Cohen--Macaulay subcurves are $B$ and $C$. The restriction sequence
\[
 0\longrightarrow\mathcal O_X(-K_X)
 \longrightarrow\mathcal O_X(K_X)
 \longrightarrow\mathcal O_C(K_X)\longrightarrow0
\]
gives $H^0(C,\mathcal O_C(K_X))=0$, since $H^0(X,K_X)=0$ and
\[
 H^1(X,-K_X)\simeq H^1(X,2K_X)^*=0.
\]
Finite duality for $B\hookrightarrow C$, followed by twisting by
$-2K_X$, yields an injection
\[
 \omega_B\otimes\mathcal O_B(-2K_X)
 \simeq\mathcal H\!om_C(\mathcal O_B,\mathcal O_C(K_X))
 \lhook\joinrel\longrightarrow\mathcal O_C(K_X).
\]
We infer that $H^0(B, \omega_B\otimes\mathcal O_B(-2K_X))=0$, and thus
$H^1(B,\mathcal O_B(2K_X))=0$ by Serre duality.  Since $B\not\subset\Bs|2K_X|$ by assumption, $H^0(X, 2K_X)\rightarrow H^0(B, \mathcal{O}_B(2K_X))$ is nonzero. Now Riemann--Roch gives
\[
 1\leq h^0(B,\mathcal O_B(2K_X))
 =2K_X\cdot B+1-p_a(B)=3-p_a(B).
\]
Thus $p_a(B)\leq2$ and $mK_X\cdot B=m\geq2p_a(B)+1$.

In every case the inequalities for
\cite[Theorem~1.1]{CFHR} hold, proving very ampleness of $mK_X|_C$.
\end{proof}
\begin{proof}[Proof of Theorem~\ref{thm: very ample}]
By Proposition~\ref{prop: local separation}, it remains to show that $|mK_X|$ separates every length-two subscheme $\xi$ supported at a point $p\in\Sigma$. By Theorem~\ref{thm: gorenstein tricanonical} and Table~\ref{tab: volume one}, one of the following holds, and in each case we apply Lemma~\ref{lem: curve trick} with the indicated $k$.
\begin{prooflist}
\item[$|3K_X|$ is base-point-free, $k=3$] If $|3K_X|$ does not separate $\xi$, then there is a reduced curve $C\in|3K_X|$ containing $\xi$ by \cite[Lemma~38$(i)$]{LR14}, and $mK_X|_C$ is very ample by Lemma~\ref{lem: tricanonical curve}.
\item[$K_X^2=\chi(X)=1$, $k=2$] Here $p_g(X)=q(X)=0$ \cite[Theorem~2.2$(ii)$]{FPR17}, so Riemann--Roch and the vanishing of $H^i(X, 2K_X)$, $i=1,2$, give $h^0(X,2K_X)=K_X^2+\chi(X)=2$. By \eqref{eq: Fujino vanishing}, $2K_X$ is globally generated at $p$. If $|2K_X|$ does not separate $\xi$, the evaluation map $H^0(X,2K_X)\to H^0(\xi,\mathcal O_\xi(2K_X))$ has rank one, and its nonzero kernel supplies a curve $C\in|2K_X|$ containing $\xi$. Since $p\in C$ and $p\notin\Bs|2K_X|$, the curve $C$ is not supported in $\Bs|2K_X|$, so $mK_X|_C$ is very ample by Lemma~\ref{lem: bicanonical curve}.
\item[$K_X^2=1$, $\chi(X)=3$] In this case $mK_X$ is very ample for every $m\geq5$ by \cite[Proposition~3.6$(v)$]{FPR17}.
\end{prooflist}
This proves the theorem.
\end{proof}

	\subsection{The pluricanonical maps of Gorenstein stable surfaces with \texorpdfstring{$K_X^2 = 1$ and $\chi(X) = 0$}{K\^{}2 = 1 and chi = 0}}\label{sec: computations}
In \cite[Theorem~A]{FPR15b} it was shown that there are exactly two Gorenstein stable surfaces with $K_X^2 = 1$ and $\chi(X) = 0$, both satisfying $p_g(X) =0$ and $q(X) = 1$. We will now explicitly compute their pluricanonical maps, thereby supplying the base point freeness of $|3K_X|$ used in the proof of   Theorem~\ref{thm: gorenstein tricanonical}.

Let us recall the description of the two surfaces $X_{0.1}$ and $X_{0.2}$ from  \cite[Theorem~A]{FPR15b}: 
	we take $\bar X=\IP^2$ with homogeneous coordinates $x,y,z$ and $\bar D=L_1+\cdots +L_4$, where 
	\[\text{$L_1=\{x=0\}$, $L_2=\{y=0\}$, $L_3=\{z=0\}$  and $L_4=\{x+y+z=0\}$.}\]
	
	Denote by $P_{(ij)}\in \bar D$ the intersection point of $L_i$ and $L_j$, explicitly,
	\begin{gather*}
		P_{(12)} = (0 :0 :1 ),\; P_{( 13 )} = (0 :1 :0 ),\; P_{( 23 )} = (1 :0 :0 ),\\
		P_{(14  )} = (0 :1 :-1 ),\; P_{(24  )} = (1 :0 :-1 ),\; P_{( 34 )} = (1 :-1 :0 ).
	\end{gather*}
	
	The normalisation of $\bar D$ is $\bar D^\nu = \bigsqcup L_i$ and  we denote by  $P_{ij}$  the point  of $L_i\subset \bar D^{\nu}$  that maps to $P_{(ij)}$, so that  each component of $\bar D^{\nu}$ contains three marked points. 

	We want to define an involution $\tau \colon \bar D^\nu\to \bar D^\nu$ such that $\tau(L_1) = L_2$ and $\tau(L_3) = L_4$ and marked points are mapped to marked points. Any such map is uniquely determined by its action on the three marked points on each component, that is, by two bijections
	\begin{gather*}
		\phi_{12}\colon \{ P_{12}, P_{13}, P_{14}\}\to \{P_{21}, P_{23}, P_{24}\},\\
		\phi_{34}\colon\{P_{31}, P_{32}, P_{34}\}\to \{P_{41},P_{42}, P_{43}\}.
	\end{gather*}
	The triple $(\bar X, \bar D, \tau)$ corresponds to a stable surface by
	Koll\'ar's gluing theorem (Section~\ref{sec: gluing}) and up to isomorphism the only two irregular Gorenstein stable surfaces with $K_X^2 = 1$ arise from this construction with the involution data given in  Table \ref{tab: 4lines}. 

	\begin{table}
					\begin{tabular}{c cc c}
				\toprule
				surface $X$ & $\phi_{12}$ & $\phi_{34}$& $\Aut(X)$ \\
				\midrule
				$X_{0.1}$
				&
				$\begin{pmatrix}
					P_{12} & P_{13} & P_{14}\\
					P_{23} & P_{21} &P_{24}
				\end{pmatrix}$
				&
				$\begin{pmatrix}
					P_{31} & P_{32} & P_{34}\\
					P_{41} & P_{43} &P_{42}
				\end{pmatrix}$
				&
				$\langle(14)(23)\rangle$
								\\
				\midrule
				$X_{0.2}$&
				$\begin{pmatrix}
					P_{12} & P_{13} & P_{14}\\
					P_{23} & P_{24} &P_{21}
				\end{pmatrix}$
				&
				$\begin{pmatrix}
					P_{31} & P_{32} & P_{34}\\
					P_{41} & P_{43} &P_{42}
				\end{pmatrix}$
				&
				$\{1\}$
				\\
				\bottomrule
			\end{tabular}
				\caption{Gluing data for Gorenstein stable surfaces with $K_X^2 = 1$ and $\chi(X) = 0 $ from  \cite{FPR15b}.}
		\label{tab: 4lines}
	\end{table}
	
		The non-normal locus of these surfaces has two components, both rational curves with three points identified to a 3-multi-node. Where the two components meet, the surfaces have a degenerate cusp, where analytically locally six components meet. Every other point is smooth or a normal crossing point. An illustration and information on their topological invariants can be found in \cite{FPR15b}. 

	Now we look at the pluricanonical systems on these surfaces.
	\begin{prop}\label{prop: chi zero surfaces}
		Let $X = X_{0.1}$ or $X = X_{0.2}$. Then
		\begin{enumerate}
			\item $|2K_X|$ has base points;
			\item $|3K_X|$ is base-point-free;
			\item $|4K_X|$ is base-point-free but not an embedding.
		\end{enumerate}
	\end{prop}
	\begin{proof}
	 By Riemann--Roch and vanishing \eqref{eq: Fujino vanishing} we have $P_m(X) = h^0(X, mK_{X}) = \binom m 2$ for $m\geq 2$. 
	 
	 \begin{prooflist}
	 \item[(i)] Since $P_2(X)=1$, every point of the unique bicanonical curve is a base point of $|2K_X|$. For both surfaces, the pullback of this curve to $\IP^2$ is the smooth conic
	 \[ x^2+xy+y^2+xz+yz+z^2=0,\]
	 that is, $\pi^*H^0(X,2K_X)$ is spanned by $(x^2+xy+y^2+xz+yz+z^2)\cdot\Omega^{\otimes2}$, with $\Omega$ as in the proof of $(ii)$ below; see \eqref{eq: Omega}. This is verified by the same computation as in $(ii)$, now with the sign condition $\tau^*\Res=\Res$. The conic meets each line $L_i$ in two points different from the marked points, so $\Bs|2K_X|$ meets the non-normal locus outside the degenerate cusp.
	 
	 \item[(iii)] $|4K_X|$ is base-point-free by \cite[Theorem~24]{LR14}, and $\dim |4K_X| = 5$ while the degenerate cusp has embedding dimension six, so $|4K_X|$ is not an embedding. 
	 
	\item[(ii)] It remains to prove $(ii)$ by computing the tricanonical map. By \eqref{eq: descent} we have
	 \begin{align*}
	 	 H^0(X, 3K_X) &\isom \pi^* H^0(X, 3K_X) \\
	 	 &= \left\{ \alpha \in H^0(\IP^2, \omega_{\IP^2}(L_1+L_2+L_3+L_4)^{\otimes 3}) \mid\tau^* \Res_{\bar D^\nu} \alpha = - \Res_{\bar D^\nu} \alpha \right\}
	 	\end{align*}
	 	 where $\tau$ is the gluing involution, $\bar D^\nu$ is the disjoint union of the four lines and $\Res$ is induced by the four residue maps $\omega_{\IP^2}(L_1+L_2+L_3+L_4)^{\otimes 3}\to \omega_{L_i}\left(\sum_j P_{ij}\right)^{\otimes 3}$.
	 	 
	 	 For clarity, let us denote by $x_i, y_i, z_i$ the restriction of the homogeneous coordinates to $L_i$, so that the corresponding coordinate rings are 
	 	 \[ B_1 = \IC[y_1, z_1], \; B_2 = \IC[x_2, z_2], \; B_3 = \IC[x_3, y_3], \; B_4 = \IC[x_4, y_4, z_4]/(x_4+y_4+z_4).\]
	 	 
	 	 With these notations we can write down the involutions $\tau_i$, that is the isomorphisms $\phi^i_{kl}$, leading to $X_{0.i}$, explicitly:
	 	 \begin{align*}
	 	 	&\tau_1^* \colon 
	 	 	&& \phi_{12}^1(y_1, z_1) = (z_2, x_2) 
	 	 	&&  \phi^1_{34}(x_3, y_3) = (x_4, -(x_4+y_4))\\
	 	 	& \tau_2^*\colon 
	 	 	&& \phi_{12}^2(y_1, z_1) = (-z_2, x_2+z_2) 
	 	 	&& \phi^2_{34}(x_3, y_3) = (x_4, -(x_4+y_4))
	 	 \end{align*}
	 	 It is straightforward to check that these maps on coordinate rings induce the bijection on marked points given in Table \ref{tab: 4lines}.

	 	 Consider now the rational 2-form on $\IP^2$ with poles along $\bar D$ given by
	 	 \begin{equation}\label{eq: Omega}
		  \Omega := \frac{ x\, dy\wedge dz - y\, dx\wedge dz + z\, dx\wedge dy}{xyz(x+y+z)}.
		  \end{equation}
	 	 Its residues along the four lines are 
	 	 \begin{align*}
	 	 	\Omega_1 & = \Res_{L_1} (\Omega) = \frac { z_1\, dy_1-y_1\, dz_1}{y_1z_1(y_1+z_1)}, 
	 	 	&
	 	 	\Omega_2 & = \Res_{L_2} (\Omega) = \frac {x_2\, dz_2 - z_2\, dx_2}{x_2z_2(x_2+z_2)},\\
	 	 	\Omega_3 & = \Res_{L_3} (\Omega) = \frac {y_3\, dx_3 -x_3\, dy_3}{x_3y_3(x_3+y_3)},
	 	 	&
	 	 	\Omega_4
	 	 	& = \Res_{L_4} (\Omega)
	 	 	= \frac {x_4\, dy_4 -y_4\, dx_4}{x_4y_4(x_4+y_4)}.
	 	 \end{align*}	 	 
	 	 A straightforward computation gives $
	 	 	\tau_1^*\Omega_1 = \Omega_2 = \tau_2^*\Omega_1$ and 
	 	 	$\tau_1^*\Omega_3  =  \Omega_4 = \tau_2^*\Omega_3$
	 	 so the residues are invariant under the involutions and we only have to consider the action on the coefficients.	 	 

With these notations in place, one can compute that, writing $\alpha = c\cdot\Omega^{\otimes 3}$ with a cubic form $c$, bases of $\pi ^* H^0(X_{0.i}, 3K_{X_{0.i}})$ are given by the following cubics 
\begin{align*}
	X_{0.1} \colon &&& 	 z\, (x-y)\, (x+y+z)
	\\
&&& 	 x\, y \,(x+y+2z) 
\\	
&&&	
 x^{3}-x\,y^{2}+y^{3}-4\,x\,y\,z+y^{2}z-x\,z^{2}-z^{3}\\
 X_{0.2} \colon &&&
z\, (x+y)\, (x+y+z)
\\
&&&
x\, y \,(x+y+2z) 
\\
&&&
x^{3}-3\,x\,y^{2}-y^{3}-5\,x\,y\,z-2\,y^{2}z-x\,z^{2}-3\,y\,z^{2}-z^{3}
\end{align*}
In both cases, the base locus of the pencil generated by the first two sections is a finite number of points, and one can check directly that the last section is non-zero at all these points. 
This concludes the proof.\qedhere
	\end{prooflist}
	\end{proof}

\section{Proofs of the main results}
We first reduce to the case without boundary and then collect the facts from the individual sections to obtain the results stated in the introduction.
\subsection{Reduction to the case without boundary}
\begin{lem}[Doubling a reduced boundary]\label{lem: double}
Let $(X,\Delta)$ be a stable log surface with $\Delta\ne0$ and global
index $I$. Then there is a stable surface $X^{\mathrm d}$ which contains $X$ as a closed subscheme, such that $IK_{X^{\mathrm d}}$ is Cartier, $\mathcal O_{X^{\mathrm d}}(IK_{X^{\mathrm d}})|_X\cong\mathcal O_X(L)$, and $K_{X^{\mathrm d}}^2=2(K_X+\Delta)^2$.
\end{lem}
\begin{proof}
Glue two copies of $(\bar X,\bar D+\bar\Delta)$ by keeping the involution $\tau$ on each copy of $\bar D^\nu$ and identifying the two copies of $\bar\Delta^\nu$ via the identity. This involution is generically fixed-point-free and preserves the different, so by Koll\'ar's gluing theorem (Section~\ref{sec: gluing}) it defines a stable surface $X^{\mathrm d}$; its boundary is empty, since every component of $\bar\Delta$ has become a conductor component. The pullback of $K_{X^{\mathrm d}}$ to each copy of $\bar X$ is $K_{\bar X}+\bar D+\bar\Delta=\pi^*(K_X+\Delta)$, which gives the volume.

On each copy of $\bar X$ the gluing makes the same identifications as for $X$. Hence the image $X_j\subset X^{\mathrm d}$ of the $j$-th copy receives a finite birational universal homeomorphism from $X$, which is an isomorphism because $X_j$ is seminormal (Section~\ref{sec: vanishing}). As $X^{\mathrm d}$ is seminormal as well, a local function on $X^{\mathrm d}$ is a pair $(f_1,f_2)$ of local functions on $X_1\cong X\cong X_2$ that agree on $\Delta$.

By \eqref{eq: descent}, a local section of $\omega_{X^{\mathrm d}}^{[I]}$ is a pair $(s_1,s_2)$ of local sections of $\mathcal O_X(L)$ with $\Res_{\bar\Delta^\nu}(s_2)=(-1)^I\Res_{\bar\Delta^\nu}(s_1)$. If $\sigma$ generates $\mathcal O_X(L)$ near a point of $\Delta$, then $\Res_{\bar\Delta^\nu}(\sigma)$ is a local generator as well, so these pairs are exactly the multiples $(f_1,f_2)\cdot(\sigma,(-1)^I\sigma)$ with $f_1=f_2$ on $\Delta$. Thus $\omega_{X^{\mathrm d}}^{[I]}$ is invertible, generated by $(\sigma,(-1)^I\sigma)$, and its restriction to $X_1\cong X$ is $\mathcal O_X(L)$.
\end{proof}
One can show that the global index of $X^{\mathrm d}$ is exactly $I$, but we will not need this.

\subsection{Proof of Theorem \ref{thm: A}}
If $\Delta=0$ and $I\geq2$, the assertion is Theorem~\ref{thm: high embedding}, and if $\Delta=0$ and $I=1$ it is Theorem~\ref{thm: very ample}.
If $\Delta\ne0$, let $X^{\mathrm d}$ be the double of Lemma~\ref{lem: double}. Its global index $I^{\mathrm d}$ divides $I$, so $mIK_{X^{\mathrm d}}=\frac{mI}{I^{\mathrm d}}\cdot I^{\mathrm d}K_{X^{\mathrm d}}$ is very ample for $m\geq5$ by the case $\Delta=0$. Restricting to $X$ shows that $mI(K_X+\Delta)$ is very ample. \qed

\subsection{Proof of Theorem \ref{thm: B}}
For $I\geq2$, this is \cite[Theorem~24]{LR14}, which needs no assumption on the volume. For $I=1$ and $\Delta=0$, use
Theorem~\ref{thm: gorenstein tricanonical} and
Proposition~\ref{prop: chi zero surfaces}. If $I=1$ and $\Delta\ne0$,
the double $X^{\mathrm d}$ in Lemma~\ref{lem: double} is connected and Gorenstein and has
volume $K_{X^{\mathrm d}}^2=2(K_X+\Delta)^2\geq2$. Thus $3K_{X^{\mathrm d}}$
is globally generated, and its sections restrict to generate $3(K_X+\Delta)$ on $X$.\qed

\section{A remark on the geography of Gorenstein stable surfaces}
In \cite{LR16} we studied the geography of Gorenstein stable surfaces and noted the inequality $P_2(X) = K_X^2+\chi(X)\geq 0$ but could not find an example where equality holds. 
	
We take this opportunity to remark that in fact there is none:
\begin{prop} Let $X$ be a Gorenstein stable surface. Then  $P_2(X)\geq 1$.
 \end{prop}
 In particular, the surfaces constructed in \cite[Section~5]{LR16} cover all possible invariants of Gorenstein stable surfaces with $\chi(X)\leq 0$.
\begin{proof}
Riemann--Roch and vanishing give
$P_2(X)=K_X^2+\chi(X)$, and $K_X^2$ is a positive integer.
Consequently there is nothing to prove when $\chi(X)\ge0$.
If $\chi(X)<0$, Blache's positivity theorem
\cite[Theorem~2]{bla94} implies that $X$ is non-normal.
For every slc centre $W$, the restriction map
\[
 r\colon H^0(X,2K_X)\longrightarrow H^0(W,2K_X|_W)
\]
is surjective by \eqref{eq: Fujino vanishing}.
If there is a zero-dimensional slc centre, its reduced structure is
a point and this already produces a nonzero bicanonical section.

Otherwise $D$ has no singular points, because these are zero-dimensional slc centres (Section~\ref{sec: vanishing}), so the conductor $D\subset X$ is a disjoint union of smooth curves.  Choose a component $D_0$ of $D$, and let
$\bar D_0\to D_0$ be its inverse image in the normalisation.
This is a finite double cover of smooth curves, possibly disconnected,
with branch divisor of degree $b_0\ge0$.  Hence adjunction gives
$\pi^*(K_X|_{D_0})=K_{\bar D_0}$, and the Hurwitz formula yields 
\[
 2d_0=4g(D_0)-4+b_0 \text{ with } 
 d_0=\deg(K_X|_{D_0})>0.
\]
Riemann--Roch on the smooth curve $D_0$ now gives
\[
 h^0(D_0,2K_X|_{D_0})\ge2d_0+1-g(D_0)
   =\frac{6d_0+b_0}{4}>0.
\]
The resulting nonzero section lifts to $X$ by the surjectivity of $r$ above. 
\end{proof}

\bibliographystyle{alpha}
\bibliography{references-5KX}

\end{document}